\documentclass[11pt,a4paper]{article}
\usepackage{enumerate, amsmath, amsthm}
\usepackage{amssymb, amsfonts, mathrsfs}

\usepackage[pdftex]{graphicx}

\def\P{{ \mathbb{P}}}

\newtheorem{thm}{Theorem}

\newtheorem{lem}[thm]{Lemma}

\newtheorem{fact}[thm]{Fact}
\newtheorem{prop}[thm]{Proposition}

\newtheorem{que}[thm]{Question}

\bfseries\normalfont

\begin{document}

\title{Spanning bipartite quadrangulations of triangulations of the projective plane}

\author{
Kenta Noguchi\thanks{Department of Information Sciences, 
Tokyo University of Science, 
2641 Yamazaki, Noda, Chiba 278-8510, Japan.
Email: {\tt noguchi@rs.tus.ac.jp}}
}

\date{}
\maketitle

\noindent
\begin{abstract}
We completely characterize the triangulations of the projective plane that admit a spanning bipartite quadrangulation subgraph. 
This is an affirmative answer to a question by K\"undgen and Ramamurthi (J Combin Theory Ser B 85, 307--337, 2002) for the projective planar case.
\end{abstract}


\section{Introduction}
Let $\Sigma$ be a \textit{surface}, that is, a compact connected $2$-dimensional manifold without boundary. 
Specifically, $\P$ denotes the projective plane. 
In this paper, the graphs may have multiple edges and basically contain no loops. 
We sometimes emphasized this fact using the term \textit{multigraphs}. 
For technical reasons, only when noted can a multigraph embedded on a surface have noncontractible loops (but no contractible loops). 
A graph with neither loops nor multiple edges is denoted as \textit{simple} graph.

The \textit{triangulation} (resp. \textit{quadrangulation}) of $\Sigma$ is a graph embedded on $\Sigma$ such that each facial closed walk has length $3$ (resp. $4$). 
When $\Sigma$ is the sphere, the graph is considered as a \textit{plane} graph. 
The coloring problem of graphs on surfaces with some constraints on faces has been widely studied in various contexts. 
For proper coloring, cyclic coloring is popular; see \cite{CHJ, NNO1} for example. 
For not necessarily proper coloring, weak (polychromatic) coloring is popular; see \cite{ABBBCSSZ, KR} for example. 
K\"{u}ndgen and Ramamurthi \cite{KR} considered the weak coloring of triangulations (or embedded graphs in general) 
and raised the following question (restated with different words). 
The \textit{weak coloring} of an embedded graph on $\Sigma$ is a (not necessarily proper) coloring of the vertices 
such that no face is monochromatic. 

\begin{que}[K\"{u}ndgen and Ramamurthi \cite{KR}, Question 11.4]
\label{mainQ}
Is there a constant $c(\Sigma)$ depending only on the surface $\Sigma$, 
such that if the edge-width of a triangulation $G$ of $\Sigma$ is at least $c$, then $G$ has a weak $2$-coloring?
\end{que}

K\"{u}ndgen and Thomassen \cite{KT}, Nakamoto, Noguchi and Ozeki \cite{NNO} independently addressed this question, 
and considered \textit{spanning quadrangulation} subgraphs of a given triangulation $G$ of $\Sigma$ 
since the following proposition holds (e.g., \cite[Proposition 7]{NNO}). 

\begin{prop}
\label{equiv}
A multitriangulation $G$ of surface $\Sigma$ has a weak $2$-coloring if and only if $G$ has a spanning bipartite quadrangulation.
\end{prop}

Therefore, their bipartiteness was investigated; 
when $G$ admits a bipartite (or nonbipartite) spanning quadrangulation $Q$? 
Although $Q$ should be bipartite when $\Sigma$ is the sphere, 
the problem becomes more difficult and interesting for nonspherical surfaces. 
For Eulerian triangulation, the following results were obtained. 
Recall that triangulation $G$ is \textit{multitriangulation} if $G$ allows multiple edges. 

\begin{thm}[K\"{u}ndgen and Thomassen \cite{KT}; cf. \cite{NNO}] 
\label{KTpp}
Let $G$ be an Eulerian multitriangulation of $\P$.
If $G$ is $3$-colorable, then every spanning quadrangulation of $G$ is bipartite.
If $G$ is not $3$-colorable, then $G$ has both a spanning bipartite quadrangulation and a spanning nonbipartite quadrangulation.
\end{thm}

\begin{thm}[K\"{u}ndgen and Thomassen \cite{KT}; cf. \cite{NNO}] 
\label{KTtorus}
Let $G$ be an Eulerian multitriangulation of the torus.
Then, $G$ has a spanning nonbipartite quadrangulation.
Furthermore, if $G$ has a sufficiently large edge-width, then $G$ has a spanning bipartite quadrangulation.
\end{thm}

Nakamoto, Noguchi and Ozeki provided a necessary and sufficient condition for Eulerian triangulation of the torus to have a spanning bipartite quadrangulation as follows: 

\begin{thm}[Nakamoto, Noguchi and Ozeki \cite{NNO}] 
\label{NNOtorus}
Let $G$ be an Eulerian multitriangulation of the torus.
Then, $G$ has a spanning bipartite quadrangulation 
if and only if 
$G$ does not have a complete graph $K_7$ as a subgraph.
\end{thm}

Since $K_7$ is uniquely embeddable on the torus up to isomorphism and the edge-width is equal to $3$ (e.g., see \cite{Ne2}), 
Theorem \ref{NNOtorus} shows that, in Theorem \ref{KTtorus}, the edge-width $4$ suffices to have a spanning bipartite quadrangulation. 
In Theorems \ref{KTpp}--\ref{NNOtorus} the assumption ``Eulerian'' is crucial. 
Considering not only Eulerian but also non-Eulerian triangulations is more involved. 
In this paper, we provide a necessary and sufficient condition for all triangulations of $\P$ to have a spanning bipartite quadrangulation, 
which gives an affirmative answer to Question \ref{mainQ} when $\Sigma = \P$. 

\begin{thm}
\label{mainthm}
Let $G$ be a simple triangulation of $\P$. 
Then, $G$ does not have a spanning bipartite quadrangulation 
if and only if 
$G$ is constructed as follows: 
Let $T = K_6$ on $\P$ and $f_1, \ldots, f_{10}$ be the faces. 
Let $T_1, \ldots, T_{10}$ be plane quasi-Eulerian triangulations with respect to a face $f_i'$ for $i\in \{1, \ldots, 10\}$ (possibly $T_i = \emptyset$ for each $i$). 
Subsequently, $G$ is constructed from $T$ and $T_i$ by pasting $f_i$ and $f_i'$ for every $i\in \{1, \ldots, 10\}$. 
\end{thm}

\begin{figure}[ht]
 \centering
 \includegraphics[width=11cm]{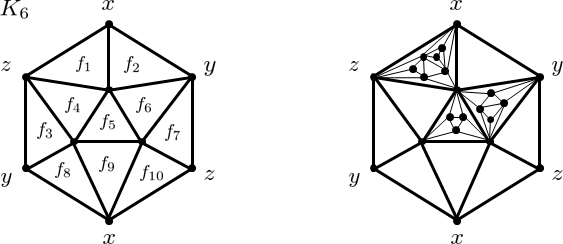}
 \caption{Left: $K_6$ on $\P$. Right: $G$ constructed from $K_6$ and $T_1; T_5; T_6$.}
 \label{fig:1}
\end{figure}

The definition of \textit{quasi-Eulerian} triangulation is somewhat complicated; see Subsection \ref{subsec2-1}. 
Such a triangulation $G$ is shown in the right-hand side of Figure \ref{fig:1} for example. 
In the figure, $T_i = \emptyset$ for $i\in \{2, 3 ,4, 7, 8, 9, 10\}$. 

Since $K_6$ is uniquely embeddable on $\P$ and the edge-width is equal to $3$ (see the left-hand side of Figure \ref{fig:1} and Fact \ref{K6}), 
Theorem \ref{mainthm} shows that, in Question \ref{mainQ}, the edge-width $c = 4$ suffices to have a weak $2$-coloring. 
Note that Theorem \ref{mainthm} is easily extended to multitriangulation of $\P$; see Section \ref{Sec2}. 

Moreover, Theorem \ref{mainthm} can be restated to the following statement in terms of the size of the bipartite subgraphs. 

\begin{thm}
\label{mainthm_restated}
Let $G$ be a triangulation of $\P$. 
Then, $G$ does not have a bipartite subgraph $H$ with $|E(H)| = \frac{2}{3}|E(G)|$ 
if and only if 
$G$ is constructed as follows: 
Let $T = K_6$ on $\P$ and $f_1, \ldots, f_{10}$ be the faces. 
Let $T_1, \ldots, T_{10}$ be plane quasi-Eulerian triangulations with respect to a face $f_i'$ for $i\in \{1, \ldots, 10\}$ (possibly $T_i = \emptyset$ for each $i$). 
Subsequently, $G$ is constructed from $T$ and $T_i$ by pasting $f_i$ and $f_i'$ for every $i\in \{1, \ldots, 10\}$. 
\end{thm}

For the lower bound of the size of bipartite subgraphs, we also prove the following theorem. 
The bound is tight by Theorem \ref{mainthm_restated}. 

\begin{thm}
\label{mainthm2}
Every triangulation $G$ of $\P$ has a bipartite subgraph $H$ with $|E(H)| \ge \frac{2}{3}|E(G)| -1$. 
\end{thm}

\section{Preliminary}
\label{Sec2}
We refer to the basic terminology in \cite{BM}, and knowledge of topological graph theory in \cite{MT}. 
Recall that a (connected) graph in which each vertex has an even degree is called \textit{Eulerian} (also called \textit{even}, for example, \cite[p.56]{BM} and \cite{MNY, NNO}). 
A \textit{$k$-vertex} is a vertex of degree $k$. 
A \textit{$k$-cycle} is a cycle of length $k$ and a multigraph can have $2$-cycles. 
Let $G$ be a ($2$-cell) embedded graph on surface $\Sigma$. 
The \textit{edge-width} of $G$ is the length of the shortest noncontractible cycle in $G$. 

A quadrangulation of $\Sigma$ is a special case of an \textit{evenly embedded graph}, 
that is, one such that each face is bounded by a closed walk of even length. 
The graph is also called \textit{even embedding}, for example, \cite{FN, No}. 
Since bipartite graphs have no cycles of odd length, any bipartite graph should be evenly embedded on any surface. 
Evenly embedded nonbipartite graphs, also called \textit{locally bipartite graphs}, have been widely studied; see \cite{LLCEHYYZ, MS}. 
The following facts are well-known for the quadrangulation of a surface. 
Note that for a quadrangulation $Q$ (or an evenly embedded graph in general) of $\P$, every noncontractible cycle of $Q$ has a length of the same parity. 

\begin{fact}
\label{bip_of_quad}
Every plane multiquadrangulation is bipartite. 
A multiquadrangulation $Q$ (or an evenly embedded multigraph in general) of $\P$ is bipartite 
if and only if 
the length of a noncontractible cycle of $Q$ is even. 
\end{fact}

Let $c$ be a $2$-coloring of $G$. 
Throughout the paper, $c$ uses two colors, black and white, $c: V(G) \to \{B, W\}$. 
(While a proper face $2$-coloring of Eulerian graphs uses two colors, red and blue.) 
Recall that $c$ is weak if no face of $G$ is \textit{monochromatic} (i.e., it receives only one color). 
Additionally, $c$ is \textit{near-weak} 
if exactly one face of $G$ is monochromatic. 
Weak $2$-coloring is also known as \textit{polychromatic $2$-coloring}. 
As Proposition \ref{equiv}, there exists a close relationship between $2$-colorings of a graph $G$ and spanning bipartite subgraphs of $G$. 
In this paper, the discussion develops while going back and forth between these two concepts. 
It is shown that Theorem \ref{mainthm2} is equivalent to the following theorem (see Section \ref{main_section} for the proof). 

\begin{thm}
\label{equiv2}
Every triangulation $G$ of $\P$ has a weak or near-weak $2$-coloring. 
\end{thm}

It is also useful to refer to the following facts and theorems. 

\begin{fact}[e.g., see \cite{Ne1}]
\label{K6}
The embedding of $K_6$ on $\P$ is unique up to isomorphism, which is a triangulation. 
In the embedding, for every edge $e$ of $K_6$, there exists a noncontractible $3$-cycle containing $e$. 
\end{fact}

\begin{thm}[K\"{u}ndgen and Ramamurthi {\cite[Corollary 4.2]{KR}}]
\label{Cor4.2}
If $K_n$ is embedded as a triangulation $G$, then $G$ has no spanning bipartite quadrangulation when $n\ge 5$.
\end{thm}

We noted the role of multiple edges. 
On a nonspherical surface, there exist two types of multiple (double) edges: a contractible $2$-cycle and noncontractible $2$-cycle. 
A contractible $2$-cycle consists of double edges that bound a $2$-cell region. 
A noncontractible $2$-cycle consists of double edges that do not bound a $2$-cell region. 
To address the operations defined in Subsection \ref{subsec2-4}, multiple edges should be allowed. 
However, some theorems and lemmas are asserted only for simple graphs to simplify the proof. 
However, most theorems and lemmas also hold for multigraphs and the proof can be modified easily. 
Theorem \ref{mainthm} for example, 
every triangulation of $\P$ with a noncontractible $2$-cycle has a spanning bipartite quadrangulation (see Lemma \ref{edge-width=2}), 
and the following proposition holds (see Subsection \ref{subsec2-1} and the proof is easy by using Lemma \ref{extend_coloring}). 

\begin{prop}
\label{contractingC}
Let $G$ be a triangulation of surface $\Sigma$ with a contractible $2$-cycle $C$. 
Contract $C$ and let $G'$ be the resulting triangulation of $\Sigma$. 
Then, $G'$ has a spanning bipartite quadrangulation if and only if $G$ has a spanning bipartite quadrangulation. 
\end{prop}

\subsection{Triangulations}
\label{subsec2-1}
In this subsection, we define several terms related to triangulations. 
Let $T$ be a triangulation of surface $\Sigma$. 
A $3$-cycle of $T$ is \textit{separating} if both the inner and outer regions of $T$ have at least one vertex. 
(Precisely, the ``inner'' region can be determined only when $T$ is drawn on the plane. 
When $\Sigma = \P$, we assume that the inner region is a $2$-cell region.) 
If $T$ has a separating $3$-cycle $C = v_1v_2v_3$, then $T$ can be divided into two smaller triangulations $T_1$ and $T_2$ 
such that $T_1$ is induced by $v_1, v_2, v_3$ and the inner vertices of $C$ and $T_2$ is induced by $v_1, v_2, v_3$ and the outer vertices of $C$. 
(Generally, the sum of the Euler genera of the two resulting surfaces equals that of $\Sigma$.) 
For the inverse operation, we say the following. 
Let $T_1$ and $T_2$ be triangulations and let $f = v_1v_2v_3 \in F(T_1)$ and $f' = v_1'v_2'v_3' \in F(T_2)$. 
Then, $T$ is said to be \textit{constructed from $T_1$ and $T_2$ by pasting $f$ and $f'$}. 
If we do not need to distinguish the three vertices when pasting, then we do not specify the labels of the vertices. 
In this case, there exist at most six possibilities for the resulting triangulations. 

If $T$ has a contractible $2$-cycle $C = v_1v_2$, 
we can contract the inner region of $C$ by deleting all inner vertices of $C$ and identifying the double edges $v_1v_2$. 
This operation is called \textit{contracting $C$}. 
When $T$ is a plane multitriangulation, by repeating the contracting $2$-cycles, we obtain the unique simple plane triangulation $T'$. 

\medskip
A simple plane triangulation $T$ with face $f$ is called a \textit{quasi-Eulerian triangulation with respect to $f$}, 
abbreviated as \textit{quasi-Eulerian triangulation w.r.t.$f$}, 
if $T$ can be constructed as follows: 
We first define the class $\mathcal{T}_1$ of simple plane triangulations. 
Next, we define $\mathcal{T}_i$ from $\mathcal{T}_{i-1}$ recursively. 
Finally, we define the class $\mathcal{T}$ of quasi-Eulerian triangulations w.r.t.$f$ by letting $\mathcal{T} = \bigcup_{i=1}^{\infty} \mathcal{T}_i$. 
Note that $\mathcal{T}_1 \subseteq \mathcal{T}_2 \subseteq \cdots$. 

\begin{figure}[ht]
 \centering
 \includegraphics[width=11cm]{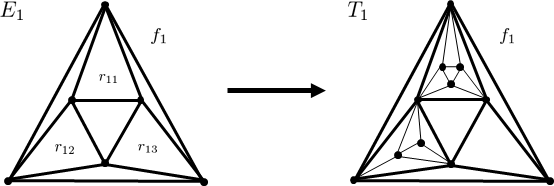}
 \caption{Quasi-Eulerian triangulation $T_1$ w.r.t.$f_1$ constructed from $E_1$ and $T_{11}; T_{12}$.}
 \label{fig:2}
\end{figure}

Let $E_1$ be a simple plane Eulerian triangulation with no separating $3$-cycles and let $f_1$ be a face of $E_1$. 
Since $E_1$ is Eulerian, $E_1$ is face $2$-colorable, using red and blue, and let $f_1$ be a red face. 
Let $r_{11}, r_{12}, \ldots, r_{1k_1}$ be all red faces of $E_1$ other than $f_1$. 
(Note that $k_1 = \frac{1}{2}|F(E_1)|-1$.) 
Let $R_{11}, \ldots, R_{1k_1}$ be plane triangulations and $r_{1j}' \in F(R_{1j})$ for $j\in \{1, \ldots, k_1\}$ (possibly $R_{1j} = \emptyset$ for each $j$). 
Subsequently, $T_1$ is constructed from $E_1$ and $R_j$ by pasting $r_{1j}$ and $r_{1j}'$ for every $j\in \{1, \ldots, k_1\}$. 
For example, see Figure \ref{fig:2}, where $f_1$ is the outer face, and $T_{13} = \emptyset$. 
Let $\mathcal{T}_1$ (w.r.t.$f_1$) be the set of all triangulations $T_1$ defined by the above procedure. 

Next, we assume that the class $\mathcal{T}_{i-1}$ is defined for some $i\ge 2$. 
Let $E_i$ be a simple plane Eulerian triangulation with no separating $3$-cycles and let $f_i$ be a face of $E_i$. 
Color the faces of $E_i$ properly by red and blue such that $f_i$ is red. 
Let $r_{i1}, r_{i2}, \ldots, r_{ik_i}$ be all red faces of $E_i$ other than $f_i$. 
Let $R_{i1}, \ldots, R_{ik_i}$ be plane triangulations and $r_{ij}' \in F(R_{ij})$ for $j\in \{1, \ldots, k_i\}$ (possibly $R_{ij} = \emptyset$ for each $j$). 
Further, 
let $b_{i1}, b_{i2}, \ldots, b_{i(k_i+1)}$ be all blue faces of $E_i$. 
Let $B_{i1}, \ldots, B_{i(k_i+1)}$ be quasi-Eulerian triangulations in $\mathcal{T}_{i-1}$ w.r.t.$b_{ij}'$ 
where $b_{ij}' \in F(B_{ij})$ for $j\in \{1, \ldots, k_i+1\}$ (possibly $B_{ij} = \emptyset$ for each $j$). 
Subsequently, $T_i$ is constructed from $E_i$ and: 
$R_{ij}$ by pasting $r_{ij}$ and $r_{ij}'$ for every $j\in \{1, \ldots, k_i\}$; 
$B_{ij}$ by pasting $b_{ij}$ and $b_{ij}'$ for every $j\in \{1, \ldots, k_i+1\}$. 
We call $E_i$ the \textit{skeleton Eulerian triangulation of $T_i$} (w.r.t.$f_i$). 
Let $\mathcal{T}_i$ (w.r.t.$f_i$) be the set of all triangulations $T_i$ defined by the above procedure. 
Note that $\mathcal{T}_i \supseteq \mathcal{T}_{i-1}$ by assuming $f_{i-1} = f_i$. 

Finally, we define $\mathcal{T}$ (w.r.t.$f$) by letting $\mathcal{T} = \bigcup_{i=1}^{\infty} \mathcal{T}_i$ (w.r.t.$f_i$), where $f \in \{f_i \mid i \ge 1 \}$. 
We omit the notation ``w.r.t.$f$'' when there is no confusion.

\subsection{Dual cubic graphs of triangulations}
\label{subsec2-2}
To investigate triangulations and quadrangulations, their dual graphs are important, and there is much research in this context. 
Let $T$ be a multitriangulation of surface $\Sigma$ and $T^*$ be its dual cubic multigraph. 
Note that an edge set $F$ of $E(T)$ corresponds naturally to an edge set $F^*$ of $E(T^*)$. 
It is well-known that $T^*$ is $2$-connected, $3$-connected if $T$ is simple, and essentially $4$-edge-connected if $T$ is simple with no separating $3$-cycles. 
Here, a $3$-connected cubic graph is \textit{essentially $4$-edge-connected} if every (edge) cut $\{e_1, e_2, e_3\}$ induces $K_{1,3}$. 
It has been shown that a separating $3$-cycle of $T$ corresponds to a $3$-cut that does not induce $K_{1,3}$. 

As explained by \cite{KR} and \cite{NNO}, 
a spanning quadrangulation $Q$ of $T$ one-to-one corresponds to a perfect matching $M^*$ of $T^*$: 
the edge set $M = E(T)-E(Q)$ is in fact $M^*$ in the dual. 
Moreover, the following theorem is presented in \cite{KR, LR, NNO} and Fact \ref{bip_of_quad}. 

\begin{thm}
\label{Thm_parity}
Let $T$ be a multitriangulation of a surface (resp. $\P$) and $M$ be an edge set of $T$ such that the edge set $M^*$ is a perfect matching in $T^*$. 
Let $Q_M = T-M$ be the spanning quadrangulation of $T$. 
Then, $Q_M$ is bipartite if and only if 
\begin{eqnarray}
\label{parity_condition}
\mbox{for any (resp. a) noncontractible cycle $C$ in $T$,} ~|E(C) \cap M| \equiv |E(C)| \pmod{2}.
\end{eqnarray}
\end{thm}

This theorem can be extended as follows: 
Let $G$ be a multigraph and $I$ be a set of nonnegative integers. 
A \textit{factor} of $G$ is a spanning subgraph of $G$ that exhibits a certain property. 
An \textit{$I$-factor} of $G$ is such that for each vertex $v$ of $G$, $d(v) \in I$. 
If $I = \{k\}$, then we call a \textit{$k$-factor} instead of a $\{k\}$-factor. 
Perfect matching is a $1$-factor. 
Let $T$ be a multitriangulation of $\Sigma$. 
Let $F^*$ be a $\{1, 3\}$-factor of $T^*$. 
Then, $T-F$ corresponds to an evenly embedded subgraph of $T$ since every face of $T-F$ is represented as a symmetric difference between an even number of triangular faces. 
Moreover, by Fact \ref{bip_of_quad}, $T-F$ is bipartite if and only if $F$ satisfies property (\ref{parity_condition}) by letting $M=F$. 
From this fact, we see that maximizing the size of a spanning bipartite subgraph of $T$ is the same as 
minimizing the size of $\{1, 3\}$-factor $M^*$ of $T^*$ with property (\ref{parity_condition}). 

Further, considering the $2$-coloring of $T$, we obtain the following fact. 

\begin{fact}
\label{Fact_factor}
Let $T$ be a triangulation $T$ of a surface $\Sigma$. 
Then, any $2$-coloring $c$ of $T$ induces a spanning bipartite subgraph $H_c = T - F_c$ of $T$, where $F_c$ is the set of monochromatic edges. 
Moreover, the dual $F_c^*$ is a $\{1, 3\}$-factor of $T^*$ with property (\ref{parity_condition}). 
Particularly, when $\Sigma$ is the sphere, 
any $\{1, 3\}$-factor $F^*$ of $T^*$ 
induces both a spanning bipartite subgraph $H = T-F$ of $T$ and 
a $2$-coloring of $T$ such that a monochromatic face $f$ corresponds to a $3$-vertex $f^*$ in $F^*$. 
\end{fact}

This can also be regarded as follows: 
$T^*-F_c^*$ is a $\{0, 2\}$-factor of $T^*$; hence, it is a union of cycles obtained by disregarding isolated vertices, 
and the regions divided by these cycles should be $2$-colorable. 
This $2$-coloring corresponds to the bipartition of the spanning bipartite subgraph $H_c$.

\subsection{Lemmas for the proof}
In this subsection, we present several lemmas. 
Lemmas \ref{monotri} and \ref{extend_coloring} were presented in our previous work \cite{NNO}. 
Lemmas \ref{extend_coloring}, \ref{extend_mono_coloring}, \ref{BWBW} and \ref{BBBW} describe the existence of a weak $2$-coloring of plane (near-)triangulations. 
A \textit{near-(multi)triangulation} is a plane (multi)graph in which each face is bounded by a $3$-cycle, except for the \textit{outer} face which is also bounded by a cycle. 

\begin{thm}[B\"{a}bler \cite{Ba}]
\label{PMe}
Let $G$ be a $2$-connected cubic multigraph. 
Then, for every edge $e$ of $G$, $G$ has a perfect matching containing $e$. 
\end{thm}

\begin{lem}[{\cite[Lemma 14]{NNO}}]
\label{monotri}
Let $G$ be an Eulerian multitriangulation of a surface.
Then, for any $2$-coloring,
$G$ has an even number of monochromatic faces. 
\end{lem} 

\begin{lem}[{\cite[Lemma 17]{NNO}}]
\label{extend_coloring}
Let $G$ be a plane multitriangulation and $f$ be a face of $G$. 
Any weak $2$-coloring of $f$ can be extended to a weak $2$-coloring of $G$. 
\end{lem}

By Lemma \ref{extend_coloring}, any weak $2$-coloring of $f$ 
can be extended to a whole weak $2$-coloring of $G$. 
To investigate when a monochromatic $2$-coloring of $f$ can be extended to a whole weak $2$-coloring, 
we used the quasi-Eulerian triangulations defined in Subsection \ref{subsec2-1}. 

\begin{lem}
\label{extend_mono_coloring}
Let $T$ be a simple plane triangulation and $f$ be a face of $T$. \\
(i) A monochromatic $2$-coloring of $f$ can be extended to a near-weak $2$-coloring of $T$ if and only if $T$ is not quasi-Eulerian with respect to $f$. \\
(ii) If $T$ is quasi-Eulerian with respect to $f$, then any monochromatic $2$-coloring of $f$ can be extended to a $2$-coloring of $T$ such that exactly two faces are monochromatic. 
\end{lem}

To prove Lemma \ref{extend_mono_coloring}, we require more lemmas. 
For a graph $G$ and $v \in V(G)$, $N_G[v]$ denotes the \textit{closed neighbor} of $v$: $N_G[v] = \{v\} \cup N_G(v)$. 

\begin{lem}
\label{K2K13factor}
Let $G = G[X, Y]$ be an essentially $4$-edge-connected bipartite cubic graph. 
Let $x\in X$ and $y\in Y$. 
Then, there exists a perfect matching in $G-(N_G[x] \cup N_G[y])$. 
\end{lem}

\begin{proof}
Let $H = G-(N_G[x] \cup N_G[y])$ and let $(X_H \subseteq X, Y_H \subseteq Y)$ be the bipartition of $H$. 
Suppose for a contradiction that $H$ does not have a perfect matching. 
By Hall's theorem (see \cite[Theorem 16.4]{BM}), 
there exists a set $S \subseteq X_H$ that satisfies $|N_H(S)| < |S|$. 

First, we suppose that $xy \not\in E(G)$. 
Let $y_1, y_2, y_3$ be the neighbors of $x$ and let $x_1, x_2, x_3$ be the neighbors of $y$. 
Consider set $N(\{x\} \cup S) = \{y_1, y_2, y_3\} \cup N(S)$. 
As $|\{y_1, y_2, y_3\} \cup N(S)| \le |\{x\} \cup S|+1$, $\{y_1, y_2, y_3\} \cup N(S)$ can have at most three edges that are joined to $X - (\{x\} \cup S)$. 
These edges form an (edge) cut that does not induce $K_{1,3}$, which is a contradiction. 

Next, suppose that $xy \in E(G)$. 
Let $y_2, y_3$ be the other neighbors of $x$ and $x_2, x_3$ be the other neighbors of $y$. 
Consider the set $A$ of edges which join $S$ and $Y - N(S)$. 
Since $xy, xy_2, xy_3, x_2y, x_3y \in E(G)$, $|A| \le 4$. 
If $|A| \le 3$, then $N(S)$ has no edge that is joined to $X - S$. 
In this case, $A$ is an (edge) cut which does not induce $K_{1,3}$, which is a contradiction. 
If $|A| = 4$, then $N(S)$ has exactly one edge $e$ that is joined to $X - S$. 
In this case, $\{e, xy\}$ is an (edge) cut, which is a contradiction. 
\end{proof}

\begin{lem}
\label{separating3}
Let $G$ be a simple plane Eulerian triangulation and $\psi: F(G) \to \{R, B\}$ be a proper face $2$-coloring of $G$. \\
(i) For any $2$-coloring of $G$, the number of monochromatic faces in $R$ equals that in $B$. \\
Moreover, if $G$ has no separating $3$-cycle, then: \\
(ii) Let $r$ and $b$ be the faces of $G$ such that $\psi(r) = R$ and $\psi(b) = B$. 
There exists a $2$-coloring of $G$ such that exactly two faces $r$ and $b$ are monochromatic.
\end{lem}

\begin{proof}
The dual $G^*$ is a $3$-connected bipartite cubic plane graph with bipartition $(R, B)$. 

(i) Let $c$ be a $2$-coloring of $G$. 
By Fact \ref{Fact_factor}, the set of monochromatic edges $F_c$ of $G$ corresponds to a $\{1, 3\}$-factor $F_c^*$ of $G^*$. 
Since $|R| = |B|$, the number of $3$-vertices of $H_c$ in $R$ and that in $B$ are the same. 
Then, the statement holds. 

(ii) Since $G$ has no separating $3$-cycle, $G^*$ is essentially $4$-edge-connected. 
Then, by Lemma \ref{K2K13factor}, $G^*$ has a $\{1, 3\}$-factor $F^*$ with exactly two $3$-vertices $r^*$ and $b^*$ and the other $1$-vertices. 
By Fact \ref{Fact_factor}, $F$ corresponds to the desired $2$-coloring of $G$ since $G$ is a plane triangulation. 
\end{proof}

A graph $G$ is \textit{bicritical} if for any two vertices $u, v$ of $G$, $G-\{u, v\}$ has a perfect matching. 
The following theorem can be proved by an argument similar to the proof of Lemma \ref{K2K13factor} 
(using Tutte's $1$-factor theorem, e.g., see \cite[Theorem 16.13]{BM}, instead of Hall's theorem). 
Therefore, we have left this proof to the reader. 
(This statement is the same as Exercise 16.4.10 in \cite[p.440]{BM}.)

\begin{thm}
\label{bicritical}
Every essentially $4$-edge-connected nonbipartite cubic graph is bicritical. 
\end{thm}

\smallskip
\begin{proof}[Proof of Lemma \ref{extend_mono_coloring}]
We prove this by induction on the number of separating $3$-cycles in $T$. 
Let $c$ be a monochromatic $2$-coloring of $f$. 

First, suppose that $T$ has no separating $3$-cycle. 
Therefore, by definition, $T$ is quasi-Eulerian w.r.t.$f$ if and only if $T$ is Eulerian. 
If $T$ is not Eulerian, then the dual $T^*$ of $T$ is essentially $4$-edge-connected nonbipartite cubic. 
Let $u^*, v^*, w^*$ be the neighbors of $f^*$ in $T^*$. 
Then, $T^*-\{u^*, v^*\}$ has a perfect matching $M$ since $T^*$ is bicritical by Theorem \ref{bicritical}. 
Let $M' = M-\{f^*w^*\} \cup K_{1,3}$ where $K_{1,3}$ consists of $f^*,u^*,v^*$ and $w^*$. 
Thus, $M'$ corresponds to a near-weak $2$-coloring of $T$ extended from $c$. 
If $T$ is Eulerian, then $c$ cannot be extended to a near weak $2$-coloring of $T$ 
since there should be a monochromatic face other than $f$ by Lemma \ref{separating3}(i). 
However, by Lemma \ref{separating3}(ii), $c$ can be extended to a desired $2$-coloring. 
In both cases, (i) and (ii) hold. 

Next, suppose that $T$ has a separating $3$-cycle $C$. 
We choose $C = v_1v_2v_3$ such that 
$T$ is constructed from $T_1$ and $T_2$ by pasting $h = v_1v_2v_3\in F(T_1)$ and $h' = v_1'v_2'v_3'\in F(T_2)$; 
$f$ is the face of $T_1$; and 
$T_2$ has separating $3$-cycles as much as possible. 

(1) If $T$ is not quasi-Eulerian w.r.t.$f$, 
then there are two cases: 
(1-1) $T_1$ is not quasi-Eulerian w.r.t.$f$. 
In this case, $c$ can be extended to a near-weak $2$-coloring $c_1$ of $T_1$ by induction. 
Since $h$ is not monochromatic in $c_1$, 
$c_1$ can be extended to a near-weak $2$-coloring of $T$ by Lemma \ref{extend_coloring}. 
(1-2) $T_1$ is quasi-Eulerian w.r.t.$f$. 
In this case, let $E$ be the skeleton Eulerian triangulation of $T_1$ with red face $f$. 
Since $T$ is not quasi-Eulerian w.r.t.$f$, $h$ should bound a blue face in $E$ and $T_2$ is not quasi-Eulerian w.r.t.$h'$. 
Since $E$ has no separating $3$-cycle, by Lemma \ref{separating3}(ii), 
there exists a $2$-coloring of $E$ such that exactly two faces $f$ and $h$ are monochromatic. 
Furthermore, since $T_2$ is not quasi-Eulerian w.r.t.$h'$, there exists a near-weak $2$-coloring of $T_2$ such that $h'$ is monochromatic by the induction hypothesis for $T_2$. 
Combining these facts and Lemma \ref{extend_coloring}, we can extend $c$ to a near-weak $2$-coloring of $T$. 

(2) If $T$ is quasi-Eulerian w.r.t.$f$, 
then we see that $T_1$ is also quasi-Eulerian w.r.t.$f$ by definition. 
Let $E$ be the skeleton Eulerian triangulation of $T_1$ with red face $f$. 
By the maximality of $T_2$, $T_2$ is attached to the face of $E$. 
For any $2$-coloring of $E$ extended from $c$, 
there exists a monochromatic blue face $b$ of $E$ by Lemma \ref{separating3}(i). 
Since the attached triangulation $B$ on $b$ is quasi-Eulerian w.r.t.$b'$ (possibly $T_2$ or empty), 
$c$ cannot be extended to a near-weak $2$-coloring of $T$ by the induction hypothesis for $B$. 
Meanwhile, $c$ can be extended to a $2$-coloring of $E$ with exactly two monochromatic faces $f$ and $b$ using Lemma \ref{separating3}(ii). 
This coloring can then be extended to a $2$-coloring of $T$ with exactly two monochromatic faces by Lemma \ref{extend_coloring} and the induction hypothesis for $B$. 

Thus, both (i) and (ii) hold by induction. 
\end{proof}

\medskip
Note that the following two lemmas can be proved without the Four Color Theorem \cite{AH} using the results of \cite{Pe} and Theorem \ref{PMe}, respectively: 
However, we used this to simplify the proof. 

\begin{lem}
\label{BWBW}
Let $G$ be a near-multitriangulation with an outer face $u_1u_2u_3u_4$. 
Then, there exists a weak $2$-coloring $c$ of $G$ such that $c(u_1) = c(u_3) = B$ and $c(u_2) = c(u_4) = W$.
\end{lem}

\begin{proof}
We added a vertex $x$ in the outer face of $G$ and join $x$ to $u_i$ for $i\in \{1, 2, 3, 4\}$. 
For the resulting plane triangulation $G'$, there exists a $4$-coloring $\phi: V(G') \to \{1, 2, 3, 4\}$ by the Four Color Theorem. 
We may assume that $\phi(x)=1, \phi(u_1)=2$ and $\phi(u_2)=3$. 
If $\phi(u_3) = 2$ (and $\phi(u_4) \in \{3, 4\}$), then let $B = \{ \phi^{-1}(\{1, 2\})\} - \{x\}$ and $W = \{ \phi^{-1}(\{3, 4\})\}$. 
If $\phi(u_3) = 4$ (and $\phi(u_4) = 3$), then let $B = \{ \phi^{-1}(\{2, 4\})\}$ and $W = \{ \phi^{-1}(\{1, 3\})\} - \{x\}$. 
This is a weak $2$-coloring of $G$ as desired. 
\end{proof}

\begin{lem}
\label{BBBW}
Let $G$ be a near-multitriangulation with an outer face $u_1u_2u_3u_4$. 
Assume that $u_1u_3 \not\in E(T)$. 
Then, there exists a weak $2$-coloring $c$ of $G$ such that $c(u_1) = c(u_2) = c(u_3) = B$ 
or a near-weak $2$-coloring $c$ of $G$ such that $c(u_1) = c(u_2) = c(u_3) = c(u_4) = B$.
\end{lem}

\begin{proof}
Let $G'$ be the multigraph obtained by identifying $u_1$ and $u_3$, and let $v (=u_1=u_3)$ be the resulting vertex. 
Since $G'$ is planar and loopless, there exists a $4$-coloring $\phi: V(G') \to \{1, 2, 3, 4\}$ by the Four Color Theorem. 
We may assume that $\phi(v) = 1$ and $\phi(u_2) = 2$. 
The $4$-coloring $\phi$ can be regarded as that of $T$ by letting $\phi(u_1) = \phi(u_3) = 1$. 
Let $B = \{ \phi^{-1}(\{1,2\})\}$ and $W = \{ \phi^{-1}(\{3, 4\})\}$. 
This is a (near-)weak $2$-coloring of $G$ as desired. 
\end{proof}

The following lemma is important for the proposed method. 

\begin{lem}
\label{existence}
Let $G$ be a multitriangulation of $\P$. 
Then, $G$ has one of the following properties: 
a $2$-vertex, $3$-vertex, $4$-vertex, $6$-vertex, or two adjacent $5$-vertices.
\end{lem}

\begin{proof}
Let $n = |V(G)|$. 
It is shown by Euler's formula that 
\begin{itemize}
\item
$G$ has a vertex of degree at most $5$, 
\item
$G$ has exactly $2n-2$ (triangular) faces, and 
\item
the following equation holds: 
\[
\sum_{v\in V(G)} (6-d_G(v)) = 6.
\]
\end{itemize}
Triangulation cannot have a vertex of degree $1$. 
We now assume that $G$ has no vertices of degree $i \in \{2, 3, 4, 6\}$. 
Let $k$ be the number of $5$-vertices. 
Then, $k \ge \frac{n}{2}$, using the above equation. 
Since any $5$-vertex is incident to five faces, the following inequality holds and there exists a face incident to at least two $5$-vertices, 
that is, $G$ has two adjacent $5$-vertices. 
\[
k\cdot 5 \geq \frac{n}{2}\cdot 5 > 2n-2
\]
\end{proof}

\subsection{Generating theorem for multitriangulations of the projective plane}
\label{subsec2-4}

\begin{figure}[ht]
 \centering
 \includegraphics[width=14cm]{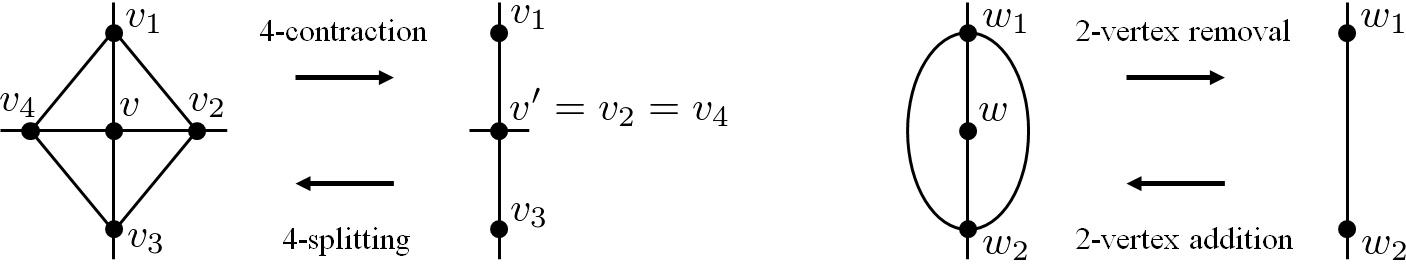}
 \caption{$4$-contraction and a $2$-vertex removal.}
 \label{ope}
\end{figure}

Let $G$ be a multitriangulation of a surface $\Sigma$ and $v$ a $4$-vertex in $G$ with neighbors $v_1, v_2, v_3, v_4$ in this order. 
Suppose that all $v_1, v_2, v_3, v_4$ are distinct. 
A \textit{$4$-contraction} of $v$ at $\{ v_2, v_4\}$ removes $v$ and identifies $v_2$ and $v_4$ and replaces two pairs of double edges with two single edges.
The inverse operation of a $4$-contraction is a \textit{$4$-splitting}.
(See the left-hand side of Figure \ref{ope}.) 
Generally, if $v_2v_4 \in E(G)$, then we do not apply the $4$-contraction of $v$ at $\{v_2, v_4\}$. 
However, in this paper we apply this when the $3$-cycle $vv_2v_4$ is noncontractible. 
Let $w$ be a $2$-vertex in $G$ and let $w_1$ and $w_2$ be the neighbors of $w$.
A \textit{$2$-vertex removal} of $w$ removes $w$ and identifies the two edges $w_1w_2$ that bound the two faces incident to $w$.
The inverse operation of a $2$-vertex removal is a \textit{$2$-vertex addition}.
(See the right-hand side of Figure \ref{ope}.)

Matsumoto et al. \cite{MNY} proved the following \textit{generating theorem} for Eulerian multitriangulations of the torus, which describes how to generate them.
Note that, by Euler's formula, every Eulerian multitriangulation $G$ of the torus has a $2$- or $4$-vertex unless $G$ is a $6$-regular triangulation. 

\begin{thm}[Matsumoto et al. \cite{MNY}] 
\label{generating}
Every Eulerian multitriangulation of the torus can be obtained from one of the $27$ graphs or a $6$-regular triangulation by a sequence of $4$-splittings and $2$-vertex additions.
\end{thm}

For $27$ graphs, see \cite[Figure 5]{MNY} or \cite[Figure 3]{NNO}. 
A generating theorem is sometimes a strong tool for investigating the properties of a family of graphs. 
In \cite{NNO}, Theorem \ref{generating} was used to prove Theorems \ref{KTtorus} and \ref{NNOtorus}. 
The existence of spanning nonbipartite quadrangulations for all minimal graphs 
guarantees the existence of those for all Eulerian multitriangulations of the torus. 
(To show the existence of spanning bipartite quadrangulations, an additional argument is required since $K_7$, one of the minimal graphs, does not have it.)
While the tool is strong, the proof depends heavily on the list of minimal graphs: 
the sporadic $27$ graphs and infinite family of $6$-regular triangulations. 
One reason for the large number of minimal graphs is that 
we sometimes cannot apply any $4$-contraction because of the small edge-width. 
In this paper we require more operations than Theorem \ref{generating}; however, there are few minimal graphs. 

Subsequently, we prepared three operations based on Lemma \ref{existence}. 
The first is a $4$-contraction; recall that in this paper, we apply it even if there exists a noncontractible $3$-cycle $vv_2v_4$. 
Note that this and the following two operations may create a noncontractible loop; however, they never create a contractible loop. 

\begin{figure}[ht]
 \centering
 \includegraphics[width=12cm]{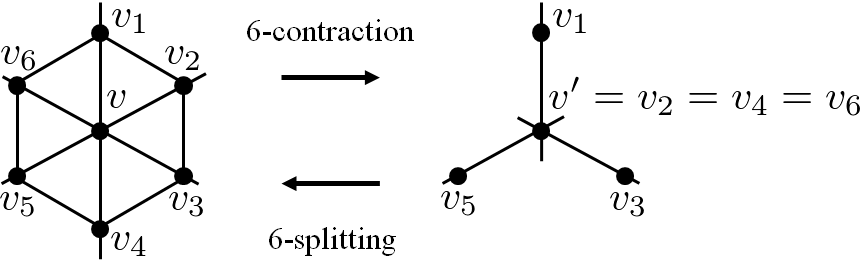}
 \caption{$6$-contraction of $v$ at $\{v_2, v_4, v_6\}$.}
 \label{ope2}
\end{figure}

Let $G$ be a multitriangulation of $\Sigma$ and let $v$ be a $6$-vertex in $G$ with neighbors $v_1, \ldots, v_6$ in this order. 
Suppose that all $v_1, \ldots, v_6$ are distinct. 
A \textit{$6$-contraction} of $v$ at $\{v_2, v_4, v_6\}$ removes $v$, identifies $v_2, v_4$ and $v_6$, 
and three pairs of double edges are replaces with three single edges.
The inverse operation of a $6$-contraction is a \textit{$6$-splitting}.
(See Figure \ref{ope2}.)
Note that if $v_2v_4 \in E(G)$ (resp. $v_2v_6 \in E(G), v_4v_6 \in E(G)$), and the $3$-cycle $vv_2v_4$ (resp. $vv_2v_6, vv_4v_6$) is contractible, 
then we do not apply the $6$-contraction of $v$ at $\{v_2, v_4, v_6\}$. 

\begin{figure}[ht]
 \centering
 \includegraphics[width=12cm]{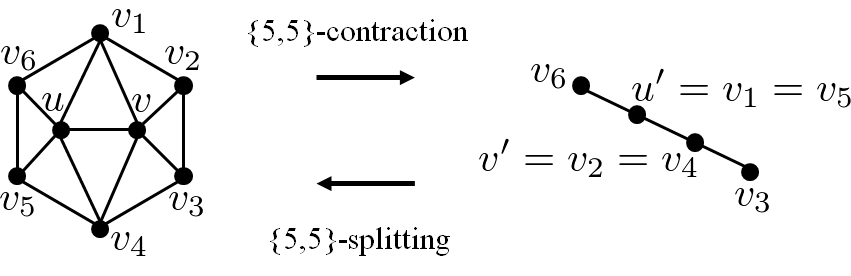}
 \caption{$\{5, 5\}$-contraction of $u; v$ at $\{v_1, v_2, v_4, v_5\}$.}
 \label{ope3}
\end{figure}

Let $G$ be a multitriangulation of $\Sigma$ 
and let $v$ and $u$ be two adjacent $5$-vertices in $G$ with neighbors $v_1, v_2, v_3, v_4$ and $v_4, v_5, v_6, v_1$ in this order, respectively. 
Suppose that all $v_1, \ldots, v_6$ are distinct, except for possibly $v_3 = v_6$. 
A \textit{$\{5, 5\}$-contraction} of $u; v$ at $\{v_1, v_2, v_4, v_5\}$ removes $u, v$, identifies $v_1$ and $v_5$, $v_2$ and $v_4$ respectively, 
and replaces three pairs of double edges with three single edges.
The inverse operation of a $\{5, 5\}$-contraction is a \textit{$\{5, 5\}$-splitting}.
(See Figure \ref{ope3}.)
Note that if $v_1v_5 \in E(G)$ (resp. $v_2v_4 \in E(G)$) and the $3$-cycle $vv_1v_5$ (resp. $vv_2v_4$) is contractible, 
then we do not apply the $\{5, 5\}$-contraction of $u; v$ at $\{v_1, v_2, v_4, v_5\}$. 
A \textit{$\{5, 5\}$-contraction} of $u; v$ at $\{v_1, v_3, v_4, v_6\}$ is similarly defined. 

The existence of spanning bipartite quadrangulations in triangulations can be preserved through three operations as described in the following lemma. 

\begin{lem} \label{reduction}
Let $G$ be a multitriangulation of a surface $\Sigma$, 
and $G'$ be a multitriangulation of $\Sigma$ obtained from $G$ by a $4$-contraction, $6$-contraction, or $\{5, 5\}$-contraction.
If $G'$ has a spanning bipartite quadrangulation, then so does $G$.
If $G'$ has a spanning nonbipartite quadrangulation, then so does $G$.
\end{lem}

\begin{proof}
Let $Q'$ be a spanning quadrangulation of $G'$ and let $M' = E(G')-E(Q')$.
We determine the edge set $M \subseteq E(G)$ such that $Q = G - M$ is a spanning quadrangulation as follows:

First, suppose that $G'$ is obtained from $G$ by a $4$-contraction of a $4$-vertex $v$. 
Let $v_1, v_2, v_3, v_4$ be the neighbors of $v$ in $G$ such that the $4$-contraction of $v$ at $\{v_2, v_4\}$ yields $G'$.
By symmetry, we may assume that one of the following holds, 
and in each case, let $M$ be the edge set as follows:
\begin{itemize}
\item $v'v_1, v'v_3\not\in M'$. Then, let $M = M'\cup \{ vv_1, vv_3\}$.
\item $v'v_1 \not\in M'$ and $v'v_3\in M'$. Then, let $M = \big( M'- \{ v'v_3\}\big) \cup \{ vv_1, v_2v_3, v_3v_4\}$.
\item $v'v_1, v'v_3\in M'$. Then, let $M = \big(M'- \{ v'v_1, v'v_3\} \big) \cup \{ v_1v_2, v_2v_3, v_3v_4, v_4v_1\}$.
\end{itemize}

Second, suppose that $G$ is obtained from $G'$ by a $6$-contraction of a $6$-vertex $v$. 
Let $v_1, \ldots, v_6$ be the neighbors of $v$ in $G$ in this order such that the $6$-contraction of $v$ at $\{ v_2, v_4, v_6\}$ yields $G'$.
By symmetry, we may assume that one of the following holds, 
and in each case, let $M$ be the edge set as follows:
\begin{itemize}
\item $v'v_1, v'v_3, v'v_5\not\in M'$. Then, let $M = M'\cup \{ vv_1, vv_3, vv_5\}$.
\item $v'v_1, v'v_3 \not\in M'$ and $v'v_5\in M'$. Then, let $M = \big( M'- \{ v'v_5\}\big) \cup \{vv_1, vv_3, v_4v_5, v_5v_6\}$.
\item $v'v_1 \not\in M'$ and $v'v_3, v'v_5\in M'$. Then, let $M = \big( M'- \{ v'v_3, v'v_5\}\big) \cup \{vv_1, v_2v_3, v_3v_4, v_4v_5, v_5v_6\}$.
\item $v'v_1, v'v_3, v'v_5\in M'$. Then, let $M = \big(M'- \{ v'v_1, v'v_3, v'v_5\} \big) \cup \{ v_1v_2, v_2v_3, v_3v_4, v_4v_5, v_5v_6, v_6v_1\}$.
\end{itemize}

Third, suppose that $G$ is obtained from $G'$ by a $\{5, 5\}$-contraction of $5$-vertices $u; v$. 
Let $v_1, v_2, v_3, v_4$ (resp. $v_4, v_5, v_6, v_1$) be the neighbors of $v$ (resp. $u$) in $G$ in this order 
such that the $\{5, 5\}$-contraction of $u; v$ at $\{ v_1, v_2, v_4, v_5\}$ yields $G'$.
By symmetry, we may assume that one of the following holds, 
and in each case, let $M$ be the edge set as follows:
\begin{itemize}
\item $v'v_3, u'v', u'v_6\not\in M'$. Then, let $M = M'\cup \{ vv_1, vv_3, uv_4, uv_6\}$.
\item $v'v_3, u'v' \not\in M'$ and $u'v_6\in M'$. Then, let $M = \big( M'- \{ u'v_6\}\big) \cup \{ vv_1, vv_3, uv_4, v_5v_6, v_6v_1\}$.
\item $v'v_3, u'v_6 \not\in M'$ and $u'v'\in M'$. Then, let $M = \big( M'- \{ u'v'\}\big) \cup \{ uv, vv_3, uv_6, v_1v_2, v_4v_5\}$.
\item $u'v' \not\in M'$ and $v'v_3, u'v_6\in M'$. Then, let $M = \big( M'- \{ v'v_3, u'v_6\}\big) \cup \{ vv_1, uv_4, v_2v_3, v_3v_4 v_5v_6, v_6v_1\}$.
\item $v'v_3 \not\in M'$ and $u'v', u'v_6\in M'$. Then, let $M = \big( M'- \{ u'v', u'v_6\}\big) \cup \{ uv, vv_3, v_1v_2, v_4v_5, v_5v_6, v_6v_1\}$.
\item $v'v_3, u'v', u'v_6\in M'$. Then, let $M = \big(M'- \{ v'v_3, u'v', u'v_6\} \big) \cup \{ uv, v_1v_2, v_2v_3, v_3v_4, v_4v_5, v_5v_6, v_6v_1\}$.
\end{itemize}

Thus, we can construct a spanning quadrangulation $Q = G - M$ in $G$ using $Q'$ in $G'$. 
Since every contraction does not change the parity of the length of the noncontractible cycle, 
it is proven by Fact \ref{bip_of_quad} that $Q$ is bipartite if and only if $Q'$ is bipartite. 
\end{proof}

\section{Proof of the main theorems}
\label{main_section}
In this section, we prove Theorems \ref{mainthm}--\ref{mainthm2} and \ref{equiv2}. 
We first show the following lemmas regarding the existence of a weak $2$-coloring of the multitriangulations of $\P$. 

\begin{lem}
\label{edge-width=1}
Let $G$ be a multitriangulation of $\P$ with a noncontractible loop $e = u_1u_1$. 
Then, there exists a weak $2$-coloring of $G$.
\end{lem}

\begin{proof}
The dual $G^*$ is $2$-connected cubic, and $G^*$ has a perfect matching $M^*$ containing $e^*$ by Theorem \ref{PMe}. 
Thus, $G-M$ corresponds to a spanning bipartite quadrangulation by Theorem \ref{Thm_parity}, and a weak $2$-coloring of $G$ by Proposition \ref{equiv}. 
\end{proof}

\begin{lem}
\label{edge-width=2}
Let $G$ be a multitriangulation of $\P$ with a noncontractible $2$-cycle $C = u_1u_2$. 
Then, there exists a weak $2$-coloring of $G$.
\end{lem}

\begin{proof}
Cut $G$ along $C$ and let $T$ be the near-multitriangulation with the outer face $C = u'_1u'_2u''_1u''_2$. 
By Lemma \ref{BWBW}, there exists a weak $2$-coloring $c$ of $T$ such that $c(u'_1) = c(u''_1) = B$ and $c(u'_2) = c(u''_2) = W$. 
Then, $c$ can be regarded as a weak $2$-coloring of $G$ by letting $c(u_1) = B$ and $c(u_2) = W$. 
\end{proof}

\begin{figure}[ht]
 \centering
 \includegraphics[width=11cm]{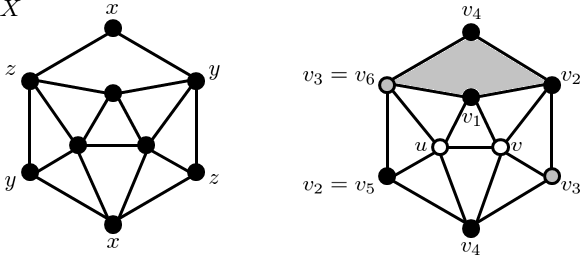}
 \caption{(Sub)graph $X$. On the right-hand side, the color of $v_3$ is either black or white.}
 \label{K6-e}
\end{figure}

\begin{lem}
\label{2-coloring_K6-e}
Let $X$ be the graph embedded on $\P$ shown in the left-hand side of Figure \ref{K6-e} ($K_6$ minus an edge). 
If a triangulation $G$ of $\P$ has a subgraph $X$ but does not have $K_6$, 
then $G$ has a weak $2$-coloring. 
\end{lem}

\begin{proof}
The right-hand side of Figure \ref{K6-e} shows that the shaded region and each triangular region may include some vertices and edges. 
Since $u, v, v_1, v_2, v_3, v_4$ form $K_6$ minus the edge $v_1v_4$, we have $v_1v_4\not\in E(G)$. 
Let $T$ be the near-triangulation as a subgraph of $G$, which consists of a shaded region with an outer face $v_1v_2v_4v_3$. 
By Lemma \ref{BBBW}, there exists a (near-)weak $2$-coloring $c$ of $T$ such that $c(v_1) = c(v_2) = c(v_4) = B$. 
Combining Lemma \ref{extend_coloring}, 
hence $c$ can be extended to a weak $2$-coloring $c'$ of $G$ as follows: 
\begin{itemize}
\item
$c'(x) = c(x)$ if $x\in V(T)$,
\item
$c'(u) = c'(v) = W$. 
\item
If a triangular region includes some vertices, then their colors by $c'$ are determined by Lemma \ref{extend_coloring}. 
\end{itemize}
Note that the color $c(v_3)$ does not affect the weakness of $c'$. 
\end{proof}

\begin{proof}[Proof of Theorem \ref{mainthm}]
First, we show the ``if'' part. 
If $G = K_6$, then $G$ is a triangulation by Fact \ref{K6} (see the left-hand side of Figure \ref{K6andX}). 
It is shown by Theorem \ref{Cor4.2} (and easy to check) that $G$ has no spanning bipartite quadrangulation.
Therefore, suppose that $G$ has $K_6$ as a proper subgraph.
Let $c$ be any $2$-coloring of $G$. 
Thus, at least one triangle of $K_6$ bounds a $2$-cell region that receives the same color (black) for its three vertices. 
The quasi-Eulerian triangulation attached to the region has another monochromatic face by Lemma \ref{extend_mono_coloring}(i). 
Hence, $c$ should have a monochromatic face of $G$. 
By Proposition \ref{equiv}, $G$ has no spanning bipartite quadrangulation. 

\begin{figure}[ht]
 \centering
 \includegraphics[width=11cm]{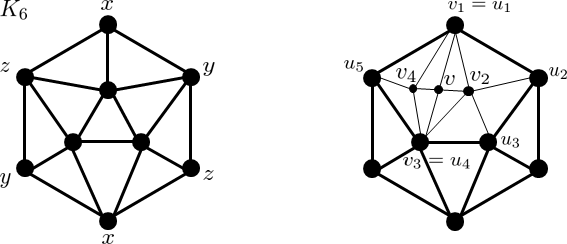}
 \caption{$K_6$ and a triangulation of $\P$.}
 \label{K6andX}
\end{figure}

Next, we show the ``only if'' part. 
First, we suppose that $G$ has $K_6$ as a subgraph. 
By the assumption, there exists at least one triangle $C$ of $K_6$ bounding a $2$-cell region 
such that $C$ and its inner vertices induce a plane triangulation that is not a quasi-Eulerian triangulation w.r.t.$C$. 
Let $c$ be a $2$-coloring of $K_6$ such that only face $C$ of $K_6$ is monochromatic; 
we color the vertices of $C$ with black and the other three vertices with white. 
By Lemmas \ref{extend_coloring} and \ref{extend_mono_coloring}(i), $c$ can be extended to a weak $2$-coloring of $G$. 
By Proposition \ref{equiv}, $G$ has a spanning bipartite quadrangulation. 

Second, suppose that $G$ does not have $K_6$ as a subgraph. 
If $G$ has a spanning bipartite quadrangulation, then we are done. 
Thus, let $G$ be a counterexample of assertion; 
that is, $G$ is a triangulation of $\P$, which does not have $K_6$ as a subgraph, 
but it has no spanning bipartite quadrangulation. 
Moreover, suppose that $G$ has the smallest number of vertices among all the counterexamples. 
There is no separating $3$-cycle in $G$ 
since, if it exists, then the removal of all inner vertices creates a smaller counterexample by Lemma \ref{extend_coloring}. 
Particularly, $G$ does not have a $3$-vertex. 
Note also that since $G$ is simple, $G$ does not have a $2$-vertex, and the edge-width of $G$ is at least $3$. 

By Lemma \ref{existence}, there exists a $4$-vertex, $6$-vertex, or two adjacent $5$-vertices in $G$.
Hence, unless Case (3-1) below, 
one of a $4$-, $6$-, or $\{5, 5\}$-contractions can be applied to $G$ to obtain a smaller multitriangulation $G'$ possibly with noncontractible loops. 
Here, if the resulting multitriangulation has contractible $2$-cycles, then contract them and make $G'$ having no contractible $2$-cycles. 
(This modification does not change the existence of spanning bipartite quadrangulation by Proposition \ref{contractingC}). 
If $G'$ does not have $K_6$ as a subgraph, then $G'$ has a spanning bipartite quadrangulation by Lemmas \ref{edge-width=1}, \ref{edge-width=2}, or the minimality of $G$. 
Hence, by Lemma \ref{reduction}, $G$ also has a spanning bipartite quadrangulation, which is a contradiction. 
Therefore, we may assume that $G'$ has $K_6$ as a subgraph. 

\medskip
Case (1) There exists a $4$-vertex $v$ in $G$. 
Let $v_1, v_2, v_3, v_4$ be the neighbors of $v$ in clockwise order. 
In this case, both $4$-contractions of $v$ at $\{v_1, v_3\}$ and $\{v_2, v_4\}$ can be applied. 
Let $G'$ (resp. $G''$) be the triangulation obtained by the $4$-contraction of $v$ at $\{v_2, v_4\}$ (resp. $\{v_1, v_3\}$) and 
$v'$ (resp. $v''$) be the resulting vertex in $G'$ (resp. $G''$). 
Subsequently, $v'$ (resp. $v''$ should be a vertex of $K_6$ in $G'$ (resp. $G''$) since $G$ does not have $K_6$. 
Let $u_1, \ldots, u_5$ be the vertices of $K_6$ other than $v'$ in $G'$ appearing clockwise around $v'$. 
By Fact \ref{K6}, $u_1, \ldots, u_5$ form a contractible $5$-cycle in $G'$ (and in $G$). 
Let $T$ be a near-triangulation with $5$-cycle $u_1u_2u_3u_4u_5$ and its inner vertices and edges in $G$. 
We may assume that one of the following three cases occurs because of the symmetry. 

(1-1) 
$\{v_1, v_3\}$ and $\{u_i, u_j\}$ coincide for some $i, j\in\{1, \ldots, 5\} ~(i \ne j)$. 
(See the right-hand side of Figure \ref{K6andX} for example.) 
In this case, the $4$-contraction of $v$ at $\{v_1, v_3\}$ creates $G''$ with the noncontractible loop $v_1v_3$ where $v_1 = v_3$, 
since there exists a noncontractible $3$-cycle $v_1vv_3$ in $G$. 
Then, $G''$ has a spanning bipartite quadrangulation by Lemma \ref{edge-width=1}, 
so does $G$ by Lemma \ref{reduction}, which is a contradiction. 

(1-2) 
$v_1$ and $u_1$ coincide, whereas $v_3$ and $u_i$ do not for any $i$. 
Since $v'u_3 \in E(G')$, $u_3v_i \in E(G)$ for some $i \in \{2, 4\}$. 
Suppose that $i=2$ (resp. $i=4$). 
We now apply the $4$-contraction of $v$ at $\{v_1, v_3\}$ in $G$. 
If $G''$ has edge-width (at most) $2$, $G''$ has a spanning bipartite quadrangulation by Lemma \ref{edge-width=2}, 
so does $G$ by Lemma \ref{reduction}, which is a contradiction. 
Therefore, $G''$ is simple (after contracting some $2$-cycles). 

Suppose that $G''$ has $K_6$. 
First, $v''$ should be a vertex of $K_6$ in $G''$ 
since $G$ does not have $K_6$. 
Note that $u_2, u_3, u_4, u_5$ are neighbors of $u_1$ in $G'$, and hence in $G$.
Second, $u_2, u_3, u_4, u_5$ should be vertices of $K_6$. 
Otherwise there exist at least two vertices $x_1, x_2$ of $K_6$ inside the contractible cycle $v''u_2u_3u_4u_5$, 
and there exists no noncontractible $3$-cycle containing edge $x_1x_2$, 
which contradicts Fact \ref{K6}. 
Third, there exists vertex $u_6$ inside the contractible cycle $u_1u_2u_3u_4u_5$ in $G$ 
such that $v'', u_2, \ldots, u_6$ form $K_6$ in $G''$. 
Specifically, $u_2u_6$ and $u_4u_6$ are edges in $G''$, and hence in $G$. 
Therefore, $u_6$ should be $v_2$ (resp. $v_4$) since the two edges $u_1v_2$ and $v_2u_3$ (resp. $u_1v_4$ and $v_4u_3$) divide $T$ into two regions. 
Therefore, $v_2$ (resp. $v_4$) should be adjacent to all $u_2, u_3, u_4, u_5$ in $G$. 
However, $u_1, \ldots, u_5$ and $v_2$ (resp. $v_4$) form $K_6$ in $G$, which is a contradiction. 

Thus, $G''$ does not have $K_6$. 
Since $|V(G'')| < |V(G)|$, 
$G''$ has a spanning bipartite quadrangulation by the minimality of $G$.
Hence, as with $G$ by Lemma \ref{reduction}, which is a contradiction.

(1-3) 
Neither $v_1$ nor $v_3$ coincides with $u_i$ for any $i$. 
In this case, since $u_1v', \ldots, u_5v' \in E(G')$, 
$v_2$ or $v_4$, say $v_2$, should be adjacent to at least three of $u_1, \ldots, u_5$. 
Without loss of generality, we may assume that $u_1v_2, u_2v_2, u_3v_2 \in E(G)$ because of the planarity of $T$. 
We now apply $4$-contraction of $v$ at $\{v_1, v_3\}$ in $G$. 
Suppose that $G''$ has $K_6$. 
Then, $u_1,\ldots,u_5$ and $v''$ should form $K_6$ in $G''$. 
However, $u_2v_1, u_2v_3\not\in E(G)$ because of the planarity of $T$, and hence $u_2v''\not\in E(G'')$, which is a contradiction. 
Thus, $G''$ does not have $K_6$ and we obtain a contradiction. 

\medskip
Case (2) There exists $6$-vertex $v$ in $G$. 
Let $v_1, \ldots, v_6$ be the neighbors of $v$ in clockwise order. 
In this case, both $6$-contractions of $v$ at $\{v_1, v_3, v_5\}$ and $\{v_2, v_4, v_6\}$ can be applied. 
Let $G'$ (resp. $G''$) be the triangulation obtained by the $6$-contractions of $v$ at $\{v_2, v_4, v_6\}$ (resp. $\{v_1, v_3, v_5\}$) 
and $v'$ (resp. $v''$) be the resulting vertex in $G'$ (resp. $G''$). 
Subsequently, $v'$ (resp. $v''$) should be a vertex of $K_6$ in $G'$ (resp. $G''$) since $G$ does not have $K_6$. 
Let $u_1, \ldots, u_5$ be the vertices of $K_6$ other than $v'$ in $G'$, 
appearing clockwise around $v'$. 
By Fact \ref{K6}, $u_1, \ldots, u_5$ form a contractible $5$-cycle in $G'$ (and hence in $G$). 
We may assume that one of the following three cases occurs because of the symmetry. 

(2-1)
At least two of $v_1, v_3, v_5$ coincide with $u_i$ and $u_j$ ($i\ne j$), say $v_1$ and $v_3$. 
In this case, 
the $6$-contraction of $v$ at $\{v_1, v_3, v_5\}$ creates $G''$ with the noncontractible loop $v_1v_3$ where $v_1 = v_3$, 
since there exists a noncontractible $3$-cycle $v_1vv_3$ in $G$. 
Subsequently, $G''$ has a spanning bipartite quadrangulation by Lemma \ref{edge-width=1}, 
and so does $G$ by Lemma \ref{reduction}, which is a contradiction. 

(2-2)
$v_1$ and $u_1$ coincide, and $v_3$ and $v_5$ do not coincide with $u_i$ for any $u_i$. 
In this case, we can proceed with the proof along with a similar argument to Case (1-2). 
Since $v'u_3 \in E(G')$, $u_3v_i \in E(G)$ for some $i \in \{2, 4, 6\}$. 
Suppose that $i=2$ (resp. $i=4, i=6$). 
We now apply the $6$-contraction of $v$ at $\{v_1, v_3, v_5\}$ in $G$. 
Thus, the same conclusion holds: 
there exists a vertex $u_6$ inside the contractible cycle $u_1u_2u_3u_4u_5$ in $G$ 
such that $v'', u_2, \ldots, u_6$ form $K_6$ in $G''$, 
and $u_6$ should be $v_2$ (resp. $v_4, v_6$). 
The only difference is that $u_1, \ldots, u_5$ and $v_i$ do not form $K_6$ in $G$ only when $i=4$. 
In this case, $G$ contains the subgraph $X$ on the left-hand side of Figure \ref{K6-e} (see Figure \ref{X}). 
Thus, $G''$ does not have $K_6$ unless $G$ contains $X$. 
By Lemma \ref{2-coloring_K6-e}, $G$ cannot contain $X$, and hence $G''$ does not have $K_6$. 
Then, $G''$ has a spanning bipartite quadrangulation by the minimality of $G$, 
and so does $G$ by Lemma \ref{reduction}, which is a contradiction.

\begin{figure}[ht]
 \centering
 \includegraphics[width=4.5cm]{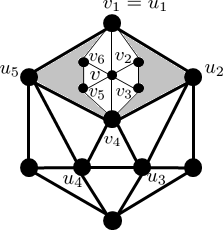}
 \caption{Subgraph $X$ in Case (2-2). The shaded regions include some edges and possibly vertices.}
 \label{X}
\end{figure}

(2-3)
Each of $v_1, v_3$ and $v_5$ do not coincide with $u_i$ for any $i$. 
If at least one of $G'$ and $G''$, say $G'$, does not have $K_6$, then $G'$ has a spanning bipartite quadrangulation by the minimality of $G$, 
and so does $G$ by Lemma \ref{reduction}, which is a contradiction. 
If both $G'$ and $G''$ have $K_6$, 
let $w_1, \ldots, w_5$ be the vertices of $K_6$ other than $v''$ in $G''$, appearing clockwise around $v''$. 
By symmetry, we may assume that $v_2, v_4$ and $v_6$ do not coincide with $w_i$ for any $i$ (otherwise this case can be reduced to Case (2-1) or (2-2)). 
Subsequently, we observed that the two sets $\{u_1, \ldots, u_5\}$ and $\{w_1, \ldots, w_5\}$ coincide. 
Now, each of $v'$ and $v''$ is adjacent to the five vertices. 
Furthermore, no $v_i$ can be adjacent to more than two $u_j$'s; 
otherwise, $v_1u_1, v_1u_2, v_1u_3 \in E(G)$ for example, $v'$ cannot be adjacent to $u_2$. 
Then, by symmetry, $G$ contains the subgraph $Y$ on the left-hand side of Figure \ref{Y}. 
Recall that every triangular region does not have any vertex or edge since there is no separating $3$-cycle in $G$. 
Let $u_3u_4v_{j+1}v_j$ be the boundary of the shaded region of $Y$. 
Now $u_3v_{j+1} \not\in E(G)$ or $u_4v_j \not\in E(G)$, and we may assume $u_3v_{j+1} \not\in E(G)$ by symmetry. 
By Lemma \ref{BBBW}, the near-triangulation induced by vertices in the shaded region with boundary $u_3u_4v_{j+1}v_j$ has a (near-)weak $2$-coloring $c$ 
such that $c(u_3) = c(u_4) = c(v_{j+1}) = B$. 
Thus, $G$ has a weak $2$-coloring, as shown in the center or right-hand side of Figure \ref{Y}. 
This contradicts that $G$ is a counterexample. 

\begin{figure}[ht]
 \centering
 \includegraphics[width=14cm]{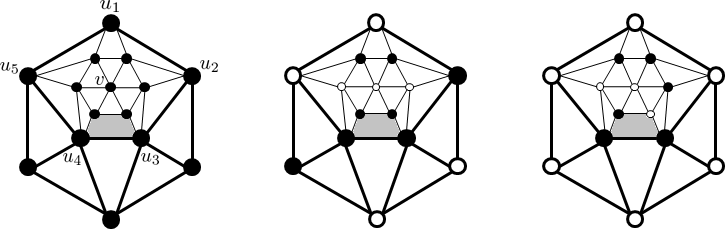}
 \caption{Subgraph $Y$ in Case (2-3). The shaded region includes some edges and possibly vertices.}
 \label{Y}
\end{figure}

\medskip
Case (3) There exist two adjacent $5$-vertices $u, v$ in $G$. 
(Assume that there are no $4$-vertices in $G$.) 
Let $v_1, v_2, v_3, v_4$ (resp. $v_4, v_5, v_6, v_1$) be the neighbors of $v$ (resp. $u$) in $G$ in this order. 
We may assume that one of the following three cases occurs because of the symmetry. 

(3-1) $v_2 = v_5$ and $v_3 = v_6$. 
In this case, $G$ has subgraph $X$; see the right-hand side of Figure \ref{K6-e}. 
Subsequently, $G$ has a weak $2$-coloring by Lemma \ref{2-coloring_K6-e}, and $G$ has a spanning bipartite quadrangulation by Proposition \ref{equiv}, which is a contradiction. 

(3-2) $v_2=v_5$ and $v_3\ne v_6$. 
In this case, $\{5, 5\}$-contraction of $u; v$ at $\{v_1, v_3, v_4, v_6\}$ can be applied. 
Let $G''$ be the resulting triangulation. 
Note that since there exist no separating $3$-cycles in $G$, 
the resulting $3$-cycle $v_2u''v''$ should be noncontractible in $G''$, and hence in $K_6$. 
By symmetry, $G$ has subgraph $Z$ or $Z'$ depicted in Figure \ref{Z}; the shaded region may include some vertices and edges. 
Furthermore, since there exists neither separating $2$- nor $3$-cycles in $G$, the shaded regions should be faces of $G$. 
Hence, $G$ is isomorphic to either $Z$ or $Z'$. 
Since both $Z$ and $Z'$ have a weak $2$-coloring, as shown in Figure \ref{Z}, 
$G$ has a spanning bipartite quadrangulation by Proposition \ref{equiv}, which is a contradiction. 
(In $Z'$, the existence of the $4$-vertex $v_6$ also contradicts the assumption.)

\begin{figure}[ht]
 \centering
 \includegraphics[width=12cm]{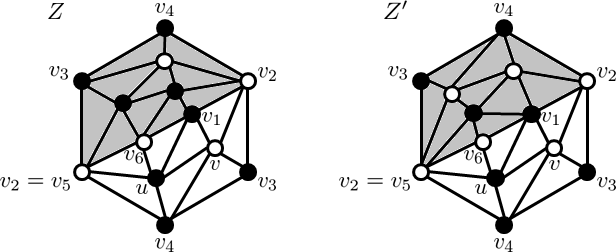}
 \caption{Graphs $Z$ and $Z'$ in Case (3-2).}
 \label{Z}
\end{figure}

(3-3) $v_2\ne v_5$ and $v_3\ne v_6$. 
In this case, we can proceed with the proof along with a similar argument to Cases (1-2) and (1-3). 
Both the $\{5, 5\}$-contractions of $u; v$ at $\{v_1, v_2, v_4, v_5\}$ and $\{v_1, v_3, v_4, v_6\}$ can be applied. 
Let $G'$ and $G''$ be the resulting triangulations. 
If at least one of $G'$ and $G''$, say $G'$, does not have $K_6$, 
then $G'$ has a spanning bipartite quadrangulation by Lemmas \ref{edge-width=1}, \ref{edge-width=2}, or the minimality of $G$, 
and so does $G$ by Lemma \ref{reduction}, which is a contradiction. 
Suppose that both $G'$ and $G''$ have $K_6$. 

(3-3a) $G$ has $K_5$ as a subgraph. 
In this case, the five vertices of $K_5$ and $v'$ (or $u'$) form $K_6$ in $G'$. 
Let $u_1, \ldots, u_5$ be the vertices of $K_6$ other than $v'$ in $G'$ appearing clockwise around $v'$. 
By Fact \ref{K6}, $u_1, \ldots, u_5$ form a contractible $5$-cycle in $G'$ (and hence in $G$). 
Note that $v_2, v_3, v_4, v_6$ may coincide with $u_i$. 
Let $T$ be a near-triangulation consisting of the $5$-cycle $u_1\cdots u_5$ and its inner vertices and edges in $G$. 
Since $u_1u', \ldots, u_5u' \in E(G')$, 
$v_1$ or $v_5$ should be adjacent to at least three of $u_1, \ldots, u_5$. 
If $v_5$, neither $u''$ nor $v''$ in $G''$ cannot be a vertex of $K_6$ in $G''$ because of the planarity of $T$. 
If $v_1$, the only possibility is that $v_1$ is adjacent to at least four of $u_1, \ldots, u_5$, say $u_1v_1, u_2v_1, u_3v_1, u_4v_1 \in E(G)$, 
and $v''(=v_1=v_3)$ is a vertex of $K_6$ in $G''$. 
(Since $G$ does not have $K_6$, $u_5v_1 \not\in E(G)$, and hence, either $u_5=v_4$ or $u_5v_5, u_5v_3 \in E(G)$ occurs; see Figure \ref{X-2}.) 
Thus, $G$ has subgraph $X$ induced by $u_1, \ldots, u_5, v_1$, which contradicts Lemma \ref{2-coloring_K6-e}. 

\begin{figure}[ht]
 \centering
 \includegraphics[width=11cm]{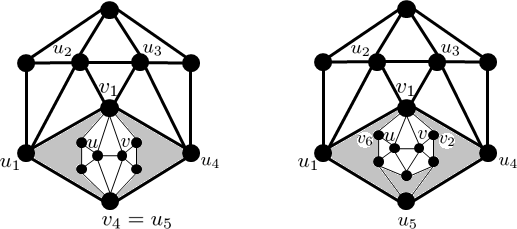}
 \caption{Subgraph $X$ in Case (3-3a). The shaded regions include some edges and possibly vertices. In the right-hand side, possibly $v_2 = u_4$ and/or $v_6 = u_1$.}
 \label{X-2}
\end{figure}

(3-3b) $G$ does not have $K_5$ as a subgraph. 
In this case, both $u'$ and $v'$ (resp. $u''$ and $v''$) should be vertices of $K_6$ in $G'$ (resp. $G''$). 
Let $u_1, u_2, u_3, u_4$ be the vertices of $K_6$ other than $u'$ and $v'$ in $G'$. 
By Fact \ref{K6}, without loss of generality, $u$ and $v$ are in the $2$-cell region bounded by a closed walk $C = u_1u_2u_3u_4u_2u_3$ in $G$; see Figure \ref{K4uv}. 
Note that $v_3$ and/or $v_6$ may coincide with $u_i$. 
Let $T$ be a near-triangulation consisting of $C$ and its inner vertices and edges in $G$ 
(to regard $T$ as a plane graph, we assume that $C$ is a $6$-cycle $u_1u_2u_3u_4u_2'u_3'$). 
Since $u_1u', u_4u' \in E(G')$, $u_1$ and $u_4$ are joined to $v_1$ or $v_5$, respectively. 
If both $u_1$ and $u_4$ are joined to the same vertex, say $v_1$, then there exist two contractible $5$-cycles $u_1u_2u_3u_4v_1$ and $u_1v_1u_4u_2u_3$ in $G$, 
one of which does not include $u$ and $v$ in the inner region; see the left-hand side of Figure \ref{K4uv}. 
%
%
That region should contain edges $u_2v_1$ and $u_3v_1$, hence, $u_1, u_2, u_3, u_4, v_1$ form $K_5$ in $G$, which is a contradiction. 
Therefore, suppose that $u_1$ and $u_4$ are joined to different vertices, that is, without loss of generality, $u_1v_1, u_4v_5 \in E(G)$; see the right-hand side of Figure \ref{K4uv}. 
By symmetry, since $u_1v', u_4v' \in E(G')$, $u_1$ and $u_4$ are joined to different vertices $v_2$ or $v_4$, respectively, and $u_1v_2, u_4v_4 \in E(G)$ because of the planarity of $T$. 

Next, consider $G''$. 
First, we show that $v_3$ and $v_6$ do not coincide with $u_i$ for any $i$. 
If $v_3 = u_1$, then there exists a separating $3$-cycle $u_1vv_1$ in $G$, which is a contradiction. 
If $v_3 = u_2$, then there exists a $4$-vertex $v_2$ bounded by a separating $4$-cycle $C' = u_1u_2vv_1$ in $G$, which contradicts the assumption. 
(Note that $C'$ contains exactly one inner vertex $v_2$ since there exist no separating $3$-cycles in $G$.) 
If $v_3 \in \{u_3, u_4\}$, then $G''$ has a noncontractible $2$-cycle $u_1u_3$ or $u_1u_4$, which contradicts Lemma \ref{edge-width=2}. 
Then, $v_3$ does not coincide with $u_i$ for any $i$, so does $v_6$ by symmetry. 
Therefore, the six vertices $u_1, u_2, u_3, u_4, u''$ and $v''$ form $K_6$ in $G''$, and it is shown by the same argument as $G'$ that $u_1v_6, u_4v_3 \in E(G)$. 
Since $u_3u' \in E(G')$ and $u_2u'' \in E(G'')$, both $u_3v_5$ and $u_2v_6$ should be edges in $G$. 
However, this is impossible because of the planarity of $T$. 

\begin{figure}[ht]
 \centering
 \includegraphics[width=11cm]{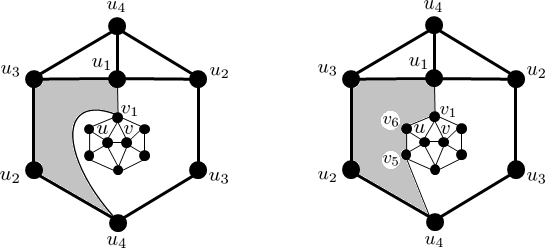}
 \caption{$K_4$ in Case (3-3b). Left: the shaded region should contain edges $u_2v_1$ and $u_3v_1$; Right: the shaded region should contain edges $u_1v_6$ and, hence, $u_3v_5$ and $u_2v_6$ but impossible.}
 \label{K4uv}
\end{figure}

From (1)--(3), such a counterexample $G$ does not exist, and the proof is completed.
\end{proof}

\begin{proof}[Proof of Theorem \ref{mainthm_restated}]
Let $G$ be a triangulation of a surface. 
Since any bipartite subgraph $H$ of $G$ cannot contain all three edges of a triangular face, 
$|E(H)| \le \frac{2}{3}|E(G)|$. 
Moreover, equality holds; that is, $G$ has a bipartite subgraph $H$ with $|E(H)| = \frac{2}{3}|E(G)|$ if and only if $G$ has a spanning bipartite quadrangulation. 
Therefore, Theorem \ref{mainthm_restated} follows from Theorem \ref{mainthm}. 
\end{proof}

\begin{proof}[Proof of Theorem \ref{equiv2}]
By Theorem \ref{mainthm_restated}, we only check that the case $G$ is constructed from $K_6$ by attaching a quasi-Eulerian triangulation to each face of $K_6$. 
Let $c$ be a near-weak $2$-coloring of $K_6$. 
By Lemmas \ref{extend_coloring} and \ref{extend_mono_coloring}(ii), $c$ can be extended to a weak or near-weak $2$-coloring of $G$. 
\end{proof}

\begin{proof}[Proof of Theorem \ref{mainthm2}]
Let $n = |V(G)|$. 
By Euler's formula, $|E(G)| = 3n-3$ and $|F(G)| = 2n-2$. 
By Fact \ref{Fact_factor}, 
if $G$ has a near-weak $2$-coloring $c$, 
then the set of monochromatic edges $F_c$ by $c$ corresponds to a $\{1, 3\}$-factor $F_c^*$ of $G^*$ that consists of exactly one $K_{1,3}$ and a matching. 
Since $|E(F^*)| = \frac{1}{2}\{3+(2n-3)\} = n$, $|E(H)| = |E(G)|-|E(F)| = 2n-3 = \frac{2}{3}|E(G)|-1$. 
Therefore, Theorem \ref{mainthm2} follows from Theorem \ref{equiv2}. 
\end{proof}

To show the equivalence of Theorems \ref{mainthm2} and \ref{equiv2}, 
finally, it is easy to show without using Theorem \ref{mainthm_restated} that Theorem \ref{mainthm2} implies Theorem \ref{equiv2}.

\section{Concluding remarks}
One of the next problems to address is the characterization of the triangulations of $\P$ that admit a spanning nonbipartite quadrangulation. 
K\"{u}ndgen and Thomassen \cite{KT} characterized them for Eulerian triangulations $G$ of $\P$: 
$G$ has a spanning nonbipartite quadrangulation if and only if $G$ is not $3$-colorable. 
We expect that the characterization for all triangulations of $\P$ is analogous to the bipartite case: 
a simple triangulation of $\P$ does not have a spanning nonbipartite quadrangulation if and only if
$G$ is constructed from a $3$-colorable triangulation of $\P$ and plane quasi-Eulerian triangulations. 
On other surfaces such as the torus or the Klein bottle, the characterization to have a spanning bipartite quadrangulation is much more involved. 
In this case, the number of nonhomotopic noncontractible curves (cycles) is at least two, and we should take care of all of them.

\section*{Acknowledgements}
Thanks to Atsuhiro Nakamoto and Kenta Ozeki for their fruitful discussions and helpful comments. 
This research was partially supported by 
JSPS KAKENHI Grant Number 21K13831.

\end{document}